\PassOptionsToPackage{reqno}{amsmath}
\documentclass[final,hidelinks]{siamart251216}
\usepackage[T1]{fontenc}
\usepackage{lmodern}
\usepackage{microtype}
\usepackage{mathtools,amssymb}
\usepackage{cite}
\hypersetup{pdftitle={Resolvent comparison for Hermitian block matrices and Kelmans transformations},pdfauthor={Xianhao Meng, Zhijun Xu, and Zhiyuan Zhang}}
\numberwithin{equation}{section}
\newsiamremark{example}{Example}
\newsiamremark{remark}{Remark}
\DeclareMathOperator{\Ran}{Ran}
\DeclareMathOperator{\rank}{rank}
\DeclareMathOperator{\diag}{diag}
\DeclareMathOperator{\adj}{adj}
\DeclareMathOperator{\spanop}{span}
\newcommand{\R}{\mathbb R}
\newcommand{\C}{\mathbb C}
\newcommand{\one}{\mathbf 1}
\newcommand{\zero}{\mathbf 0}
\newcommand{\qt}{\widetilde q}
\newcommand{\st}{\widetilde s}
\newcommand{\bt}{\widetilde\beta}
\newcommand{\U}{\mathcal U}
\newcommand{\K}{\mathcal K}
\newcommand{\leew}{\leq_{\mathrm{ew}}}
\newcommand{\geew}{\geq_{\mathrm{ew}}}
\newcommand{\colvec}[2]{\begin{bmatrix}#1\\#2\end{bmatrix}}

\makeatletter
\def\@doauth#1\and{\let\@nextdoa=\@doauth
  \def\@tempc{#1}\futurelet\@tempa\@doautha}
\makeatother
\headers{Resolvent comparison for Hermitian blocks}{Xianhao Meng, Zhijun Xu, and Zhiyuan Zhang}
\title{Resolvent comparison for Hermitian block matrices and Kelmans transformations\thanks{\textbf{Funding:} The research of all authors was supported by the Fundamental and Interdisciplinary Disciplines Breakthrough Plan of the Ministry of Education of China (JYB2025XDXM207). The third author was also supported by the National Science Foundation for Young Students of China (PhD candidate) (No.~125B2018).}}
\author{Xianhao Meng\thanks{School of Mathematical Sciences and LPMC, Nankai University, Tianjin 300071, China. China. Email: \mbox{\email{xianhaomeng@163.com}}.}
\and Zhijun Xu\thanks{School of Mathematical Sciences and LPMC, Nankai University, Tianjin 300071, China. Email: \mbox{\href{mailto:zhj_xu@mail.nankai.edu.cn}{zhj\_xu@mail.nankai.edu.cn}}.}
\and Zhiyuan Zhang\thanks{Center for Combinatorics and LPMC, Nankai University, Tianjin 300071, China. Email: \mbox{\email{jsam0331@163.com}}.}}
\begin{document}
\maketitle
\begin{abstract}
A Kelmans transformation increases or preserves the largest eigenvalue of a symmetric nonnegative matrix. We determine when this comparison holds for Hermitian block extensions with a fixed coupling and an arbitrary common complementary block. An orthogonal change of basis and Schur complementation give an exact formula for the compressed resolvent difference as a difference of two rank-one matrices. This formula characterizes admissible coupling subspaces by proportionality of projected resolvent vectors. Polynomial coefficient kernels describe all such subspaces, and cyclic subspaces give dimension bounds and a maximum-dimension result for a two-block family. We also prove a general block comparison theorem and give monotonicity and equality criteria for nonnegative rank-one couplings. Applications concern graph joins, complex attachments, and Hermitian adjacency matrices.
\end{abstract}
\section{Introduction}

The largest eigenvalue of a Hermitian matrix is monotone in the Loewner order, and the spectral radius of a nonnegative matrix is monotone in the entrywise order~\cite{HornJohnson,BermanPlemmons}. For a Hermitian block matrix, the effect of replacing a principal block also depends on its coupling to the remaining coordinates. When the coupling has signed or complex entries, both the eigenvalues and the eigenspaces of the principal block enter the comparison.

Low-rank couplings occur in several matrix constructions. Fiedler's lemma~\cite{Fiedler} and its extensions~\cite{Saravanan} determine eigenvalues of coupled blocks, and joins of normal matrices with constant row sums are treated in~\cite{Minac}. For graph matrices, spectra of joins and attachments can be obtained using coronals~\cite{McLemanMcNicholas}, vertex-subset attachments~\cite{RajkumarGayathri}, and walk matrices~\cite{Cardoso}; related formulas for deformed Laplacians appear in~\cite{Diaz}. Here we ask when replacing one block preserves the largest-eigenvalue inequality for a fixed coupling space and every choice of the complementary matrix.

We use the Kelmans transformation on symmetric nonnegative matrices, defined in Section~\ref{sec:kelmans} by replacing paired entries by their maximum and minimum. This operation extends the graph transformation~\cite{Kelmans}; related matrix transformations are studied in~\cite{KaoWeng,ChengWeng}. On the adjacency matrix of a simple graph it gives the usual graph Kelmans transformation, whose spectral-radius monotonicity is proved in~\cite{Csikvari09,Csikvari11}. We study the comparison when the transformed matrix is a principal block of a Hermitian matrix.

Our main question is: \emph{which fixed complex coupling subspaces preserve the largest-eigenvalue inequality under Kelmans transformations for every common Hermitian complementary block, and how large can such subspaces be?}

Schur complements are well suited to the comparison of Hermitian block extensions. They reduce it to the Loewner order of compressed resolvents. Their positivity and inertia properties are collected in~\cite{HornZhang}; Li and Mathias~\cite{LiMathias} obtain eigenvalue inequalities from extremal characterizations of the generalized Schur complement, and Tretter~\cite{Tretter} develops the spectral theory of block operators. Related elimination constructions occur in Kron reduction~\cite{DorflerBullo}, rational isospectral reduction~\cite{BunimovichWebb}, and reductions for restricted quantum walks~\cite{KemptonTolbert}. Entrywise bounds for inverses of $H$-matrices in terms of comparison $M$-matrices are studied in~\cite{Neumaier}. For signed or complex couplings, the relevant condition is positive semidefiniteness of the compressed resolvent difference.

The dimension question concerns subspaces on which the resolvent difference is positive semidefinite for every spectral parameter in a half-line. Range inclusion and factorization~\cite{Douglas} govern rank-one order, while higher-rank numerical ranges and totally isotropic subspaces~\cite{LiSze} concern the dimensions of constrained compressions. Rational vector spaces~\cite{Forney} and polynomial methods for matrix pencils~\cite{NoferiniVanDooren} offer coefficient descriptions of such relations. Cyclic subspaces and walk matrices~\cite{Godsil} describe the spans of resolvent vectors, and univariate sums of squares~\cite{MarshallSOS} give finite certificates for scalar nonnegativity.

The main results are in Sections~\ref{sec:kelmans} and~\ref{sec:subspaces}. An orthogonal change of basis gives two bordered matrices with a common principal block; the transformation replaces the border vector by its entrywise absolute value. The Schur-complement formula then expresses the compressed resolvent difference as a difference of two rank-one matrices (Theorem~\ref{thm:increment}). For a proper Kelmans transformation, this difference is entrywise nonnegative and has one positive and one negative eigenvalue. Its compression is positive semidefinite precisely when the two projected vectors are proportional with a suitable bound on the factor. Clearing denominators gives a description of all admissible coupling subspaces by polynomial coefficient kernels and a scalar polynomial inequality (Theorem~\ref{thm:subspace-classification}). When the common block has two diagonal blocks joined by a rank-one matrix, the resolvent equations give a dimension bound in terms of two cyclic subspaces. Let their dimensions be $k_T$ and $k_c$, respectively. Theorem~\ref{thm:rank-one-connection} determines the maximum admissible dimension when $k_T\geq k_c$ and the unique subspace attaining it when $k_T>k_c$.

Both results use a comparison theorem for general Hermitian block matrices. Theorem~\ref{thm:universal} characterizes the largest-eigenvalue inequality for every common complementary block by two conditions: comparison of the largest eigenvalues of the internal blocks, and the Loewner order of their compressed resolvents. If the resolvent difference has a negative eigenvalue, we construct an open ball of complementary blocks giving the reverse strict inequality and determine its largest possible radius. For nonnegative rank-one couplings, the required resolvent inequality follows from entrywise comparison. Generalized Schur complements and rank-one downdates then give equality criteria.

The applications include rank-one graph joins and connected graphs admitting arbitrarily many independent complex coupling vectors. We also apply the comparison to symmetric modules in Hermitian adjacency matrices~\cite{GuoMohar,LiuLi,Mohar}, both in a digraph and in its arc complement.

The paper is organized as follows. Section~\ref{sec:preliminaries} collects notation and matrix preliminaries. Section~\ref{sec:kelmans} derives the resolvent formula for the Kelmans transformation and its strictness properties. Section~\ref{sec:block-comparison} characterizes comparison of general Hermitian block extensions. Section~\ref{sec:rank-one} establishes monotonicity and equality criteria for nonnegative rank-one couplings. Section~\ref{sec:subspaces} classifies coupling subspaces and determines optimal dimensions for the two-block family. Section~\ref{sec:applications} gives the graph applications.
\section{Preliminaries}\label{sec:preliminaries}

\subsection{Notation and Hermitian eigenvalues}
For a matrix $H$, $H^T$ and $H^*$ denote its transpose and conjugate transpose, respectively. A nonnegative matrix is real and entrywise nonnegative. For real matrices of the same size, or real vectors of the same dimension, $X\leew Y$ means entrywise inequality. For Hermitian matrices of the same order, $X\preceq Y$ means that $Y-X$ is positive semidefinite, whereas $X\prec Y$ means that it is positive definite. Absolute values of matrices and vectors are entrywise.

Eigenvalues of a Hermitian matrix $H$ of order $n$ are arranged as $\lambda_1(H)\geq\cdots\geq\lambda_n(H)$, and $\lambda_{\max}(H)=\lambda_1(H)$. Its spectral radius is $\rho(H)$, and $n_+(H),n_-(H)$ count its positive and negative eigenvalues with multiplicity. The norm is Euclidean for vectors and the induced operator norm for matrices. The symbols $H^\dagger$, $\Ran H$, and $H[\mathcal T]$ denote the Moore--Penrose inverse, column space, and principal submatrix indexed by $\mathcal T$, respectively. The symbols $I$, $\one$, and $\zero$ denote the identity matrix, the all-ones vector, and a zero matrix or vector, with sizes determined by context. Matrices and explicit row or column vectors are written in square brackets. Orthogonality in $\C^n$ is with respect to $y^*z$; in particular, $w^\perp=\{z:w^*z=0\}$.

\begin{samepage}
The min--max principle gives the following inequalities for Hermitian matrices $H,K$ of order $n$ and $1\leq j\leq n$~\cite{HornJohnson}:
\begin{equation}\label{eq:variational-facts}
 H\preceq K\ \Longrightarrow\ \lambda_j(H)\leq\lambda_j(K),
 \qquad |\lambda_j(H)-\lambda_j(K)|\leq\|H-K\|.
\end{equation}
\end{samepage}
If $V\in\C^{n\times k}$ satisfies $V^*V=I$, then the interlacing lemma gives
\begin{equation}\label{eq:interlacing}
 \lambda_j(H)\geq\lambda_j(V^*HV)\geq\lambda_{j+n-k}(H),
 \qquad 1\leq j\leq k.
\end{equation}

\subsection{Schur complements and rank-one order}
The Schur-complement facts below are standard~\cite{HornZhang,LiMathias}. We include the singular case, which will be used in the equality criteria.

\begin{samepage}
\begin{lemma}\label{lem:schur}
Let $K,D$ be Hermitian matrices, with $K\succ\zero$, and $X$ have the dimensions required in
\[
 G=\begin{bmatrix}K&-X\\-X^*&D\end{bmatrix},\qquad S=D-X^*K^{-1}X.
\]
Then $G\succeq\zero$ if and only if $S\succeq\zero$, and $G\succ\zero$ if and only if $S\succ\zero$. Moreover, $n_-(G)=n_-(S)$, $\dim\ker G=\dim\ker S$, and, when $S$ is invertible, the leading block of $G^{-1}$ is
\begin{equation}\label{eq:block-inverse}
 K^{-1}+K^{-1}XS^{-1}X^*K^{-1}.
\end{equation}
For an arbitrary Hermitian matrix $K$, including a singular one,
\begin{equation}\label{eq:generalized-schur}
 G\succeq\zero\quad\Longleftrightarrow\quad
 K\succeq\zero,\quad (I-KK^\dagger)X=\zero,\quad D-X^*K^\dagger X\succeq\zero.
\end{equation}
\end{lemma}
\end{samepage}
\begin{proof}
When $K\succ\zero$, the factorization
\[
 G=\begin{bmatrix}I&\zero\\-X^*K^{-1}&I\end{bmatrix}
 \begin{bmatrix}K&\zero\\\zero&S\end{bmatrix}
 \begin{bmatrix}I&-K^{-1}X\\\zero&I\end{bmatrix}
\]
gives the positivity and inertia assertions by congruence, and its inverse gives~\eqref{eq:block-inverse}. If $G\succeq\zero$, then $K\succeq\zero$. For $y\in\ker K$ and any $z$, this quadratic form at $\colvec{ty}z$ equals $z^*Dz-2\Re(\overline t\,y^*Xz)$. Its nonnegativity for every $t\in\C$ implies $y^*Xz=0$. Thus $\Ran X\subseteq\Ran K$, which is the second condition in~\eqref{eq:generalized-schur}. Under this range condition, let $S_\dagger=D-X^*K^\dagger X$. Then
\[
 G=\begin{bmatrix}I&\zero\\-X^*K^\dagger&I\end{bmatrix}
 \begin{bmatrix}K&\zero\\\zero&S_\dagger\end{bmatrix}
 \begin{bmatrix}I&-K^\dagger X\\\zero&I\end{bmatrix}.
\]
The outer factors are invertible and are adjoints of each other. Hence $G\succeq\zero$ if and only if $K\succeq\zero$ and $S_\dagger\succeq\zero$. The same factorization proves sufficiency under the three conditions in~\eqref{eq:generalized-schur}.
\end{proof}

We shall also use the following criteria for positive semidefiniteness of a difference of outer products and a rank-one downdate.

\begin{samepage}
\begin{lemma}\label{lem:rank-one-order}
For $a,b\in\C^m$, the inequality $aa^*-bb^*\succeq\zero$ holds if and only if $b=\theta a$ for some $\theta\in\C$ with $|\theta|\leq1$. In addition, for $S\succeq\zero$ and $\eta>0$,
\begin{equation}\label{eq:downdate}
 S-\eta bb^*\succeq\zero
 \quad\Longleftrightarrow\quad
 b\in\Ran S\ \text{ and }\ \eta b^*S^\dagger b\leq1.
\end{equation}
\end{lemma}
\end{samepage}
\begin{proof}
If $aa^*-bb^*\succeq\zero$, its quadratic form on $a^\perp$ shows that $b\in\spanop\{a\}$. If $a=\zero$, then $b=\zero$; otherwise substitution gives $1-|\theta|^2\geq0$. These conditions also imply the inequality. For~\eqref{eq:downdate}, apply~\eqref{eq:generalized-schur} to $\left[\begin{smallmatrix}S&-\sqrt\eta\,b\\-\sqrt\eta\,b^*&1\end{smallmatrix}\right]$, and also take the Schur complement of its lower-right $1\times1$ block.
\end{proof}

\subsection{Nonnegative resolvents and cyclic subspaces}
The support graph of a symmetric nonnegative matrix $F$ has an edge $ij$ whenever $i\ne j$ and $F_{ij}>0$. Such a matrix is irreducible precisely when this graph is connected. The Perron--Frobenius theorem gives $\lambda_{\max}(F)=\rho(F)$; if $F$ is irreducible, an eigenvector corresponding to the largest eigenvalue can be chosen strictly positive and that eigenvalue is simple~\cite{BermanPlemmons}.

For $x>\rho(F)$, the Neumann series also identifies which entries of the resolvent are strictly positive.

\begin{samepage}
\begin{lemma}\label{lem:nonnegative-resolvent}
If $F$ is symmetric and nonnegative, then for $x>\rho(F)$,
\[
 (xI-F)^{-1}=\sum_{k\geq0}\frac{1}{x^{k+1}}F^k\geew\zero.
\]
Its $(i,j)$ entry is positive if and only if $i,j$ belong to the same component of the support graph of $F$.
\end{lemma}
\end{samepage}
\begin{proof}
The Neumann series converges because $\rho(F)/x<1$, and all its terms are entrywise nonnegative. For $i\ne j$, a positive term in $(F^k)_{ij}$ exists exactly when there is a walk from $i$ to $j$ using positive entries. A path in the support graph supplies such a term. For $i=j$, $F^0=I$, so the term $\frac{1}{x}I$ gives a strictly positive diagonal entry.
\end{proof}

For a Hermitian matrix $T$ of order $r$ and a vector $a\in\C^r$, its cyclic subspace generated by $a$ is
\[
 \K(T,a)=\spanop\{a,Ta,\ldots,T^{r-1}a\}.
\]
It is invariant under $T$ by Cayley--Hamilton; since $T$ is Hermitian, its orthogonal complement is invariant as well. The matrix with these spanning columns is a Krylov matrix. For a graph adjacency matrix and a vertex-subset characteristic vector, it is the corresponding walk matrix~\cite{Godsil}.

The following lemma relates this finite-dimensional span to resolvent vectors and adjugate coefficients.

\begin{samepage}
\begin{lemma}\label{lem:resolvent-span}
If $J$ is a nonempty open real interval disjoint from the spectrum of $T$, then
\[
 \spanop\{(xI-T)^{-1}a:x\in J\}=\K(T,a).
\]
This space also equals the span of the coefficients of the vector polynomial $\adj(xI-T)a$.
\end{lemma}
\end{samepage}
\begin{proof}
Invariance implies that each resolvent vector belongs to $\K(T,a)$. If a vector $z$ is orthogonal to these vectors on $J$, the rational function $z^*(xI-T)^{-1}a$ is zero on $J$. Multiplication by $\det(xI-T)$ gives a polynomial with infinitely many zeros, so this rational function is identically zero. Its convergent expansion $\sum_{k\geq0}x^{-k-1}z^*T^ka$ for $|x|>\|T\|$ then gives $z^*T^ka=0$ for every $k\geq0$. Taking orthogonal complements proves the reverse inclusion. Moreover, $z^*\adj(xI-T)a=\det(xI-T)z^*(xI-T)^{-1}a$ is identically zero if and only if the rational function is identically zero. This proves the assertion about coefficient spans.
\end{proof}
\begin{samepage}
\section{Kelmans transformations and compressed resolvents}\label{sec:kelmans}

\subsection{The matrix transformation in an orthonormal basis}
Let $A=[a_{ij}]\in\R^{n\times n}$ be symmetric and entrywise nonnegative, and fix distinct indices $u,v$. Let $W=\{1,\ldots,n\}\setminus\{u,v\}$. The \emph{Kelmans transformation of $A$ toward $u$ from $v$} is the symmetric matrix $\widetilde A$ defined by
\end{samepage}
\begin{align*}
 \widetilde a_{uu}&=\max\{a_{uu},a_{vv}\},&\widetilde a_{vv}&=\min\{a_{uu},a_{vv}\},&\widetilde a_{uv}&=a_{uv},\\
 \widetilde a_{uw}&=\max\{a_{uw},a_{vw}\},&\widetilde a_{vw}&=\min\{a_{uw},a_{vw}\}\qquad(w\in W),
\end{align*}
with $\widetilde A[W]=A[W]$. The Kelmans transformation is \emph{proper} if the difference vector
\[
 [a_{uu}-a_{vv},\ (a_{uw}-a_{vw})_{w\in W}]
\]
has both a positive and a negative entry. Otherwise $\widetilde A=A$ or $\widetilde A$ is obtained from $A$ by interchanging $u$ and $v$.

Throughout the paper, Kelmans transformation refers to this matrix operation. For a simple graph $G$, applying it to $A(G)$ gives the adjacency matrix of a graph $\widetilde G$ on the same vertex set. Equivalently, replace each edge $vw$ by $uw$ for $w\in N_G(v)\setminus(N_G(u)\cup\{u\})$, retaining all other edges, including $uv$ when present. Here $N_G(x)$ denotes the neighborhood of $x$ in $G$. This is the graph Kelmans transformation toward $u$ from $v$~\cite{Csikvari11}. It is proper precisely when both $N_G(u)\setminus(N_G(v)\cup\{v\})$ and $N_G(v)\setminus(N_G(u)\cup\{u\})$ are nonempty.

Let $d=n-1$ and use the orthonormal vectors
\[
 e_+=\frac{1}{\sqrt2}(e_u+e_v),\qquad e_-=\frac{1}{\sqrt2}(e_u-e_v),\qquad
 P=[e_+,\ (e_w)_{w\in W}],\quad Q=[P,e_-].
\]
The columns of $P$ form an orthonormal basis of $(e_u-e_v)^\perp$. Order the columns of $P$ as follows: column $1$ is $e_+$, and column $j$ is $e_{w_j}$ for $2\leq j\leq d$, where $W=\{w_2,\ldots,w_d\}$ with $w_2<\cdots<w_d$. In this basis,
\begin{equation}\label{eq:normal-form}
 Q^TAQ=\begin{bmatrix}F&q\\q^T&\delta\end{bmatrix},\qquad
 Q^T\widetilde AQ=\begin{bmatrix}F&\qt\\\qt^T&\delta\end{bmatrix},\qquad \qt=|q|,
\end{equation}
where
\begin{align*}
 F_{11}&=\frac{a_{uu}+a_{vv}}2+a_{uv},&F_{1j}&=\frac{a_{u w_j}+a_{v w_j}}{\sqrt2},&F[\{2,\ldots,d\}]&=A[W],\\
 q_1&=\frac{a_{uu}-a_{vv}}2,&q_j&=\frac{a_{u w_j}-a_{v w_j}}{\sqrt2},&\delta&=\frac{a_{uu}+a_{vv}}2-a_{uv}.
\end{align*}
Here $2\leq j\leq d$, and the rows and columns $2,\ldots,d$ of $F$ correspond to the rows and columns $w_2,\ldots,w_d$ of $A$, in that order. The block $F$ is nonnegative and unchanged by the transformation. Properness is equivalent to $q$ having both signs.

A common scalar shift makes the diagonal entries of both bordered matrices nonnegative. Entrywise comparison then gives the following bounds.

\begin{samepage}
\begin{lemma}\label{lem:inverse-comparison}
The spectral radii satisfy $\rho(\widetilde A)\geq\rho(A)$. For $x>\rho(\widetilde A)$,
\[
 \bigl|Q^T(xI-A)^{-1}Q\bigr|
 \leew Q^T(xI-\widetilde A)^{-1}Q.
\]
\end{lemma}
\end{samepage}
\begin{proof}
Denote the two matrices in~\eqref{eq:normal-form} by $B,\widetilde B$. Choose $t\geq0$ large enough that both $B+tI$ and $\widetilde B+tI$ are positive semidefinite. Their diagonal entries are then nonnegative, so $|B+tI|\leew \widetilde B+tI$. For every real unit vector $y$,
\[
 y^T(B+tI)y\leq |y|^T(\widetilde B+tI)|y|.
\]
Taking the maximum over $y$ gives $\lambda_{\max}(B)\leq\lambda_{\max}(\widetilde B)$. Since $A,\widetilde A$ are nonnegative and symmetric, this is the asserted spectral-radius comparison. If $x>\rho(\widetilde A)$, then $x+t$ exceeds the spectral radii of both shifted matrices. Entrywise comparison of their powers gives
\[
 \left|\sum_{k\geq0}\frac{1}{(x+t)^{k+1}}(B+tI)^k\right|
 \leew\sum_{k\geq0}\frac{1}{(x+t)^{k+1}}(\widetilde B+tI)^k.
\]
The two convergent series are $(xI-B)^{-1}$ and $(xI-\widetilde B)^{-1}$, respectively, which proves the asserted inverse inequality.
\end{proof}

\subsection{The exact difference}
For $x>\rho(\widetilde A)$, define
\begin{equation}\label{eq:schur-data}
 M=(xI-F)^{-1},\qquad s=x-\delta-q^TMq,\qquad
 \st=x-\delta-\qt^TM\qt.
\end{equation}
All three quantities depend on $x$; this argument is omitted from the notation within formulas. Since $F=P^TAP$, the interlacing lemma~\eqref{eq:interlacing} and Lemma~\ref{lem:inverse-comparison} give $\rho(F)\leq\rho(A)\leq\rho(\widetilde A)<x$. Thus $xI-F\succ\zero$, and Lemma~\ref{lem:nonnegative-resolvent} gives $M\geew\zero$. Applying Lemma~\ref{lem:schur} to $xI-Q^TAQ$ and $xI-Q^T\widetilde AQ$ gives $s,\st>0$. Moreover, $M\geew\zero$ and $\qt=|q|$ imply $q^TMq\leq\qt^TM\qt$, hence $0<\st\leq s$.
Let
\[
 R_A(x)=P^T(xI-A)^{-1}P,\qquad \Delta(x)=R_{\widetilde A}(x)-R_A(x).
\]

\begin{samepage}
\begin{theorem}\label{thm:increment}
For every $x>\rho(\widetilde A)$,
\begin{equation}\label{eq:rank-two}
 R_A=M+\frac{1}{s}(Mq)(Mq)^T,\qquad
 \Delta=pp^T-rr^T,\qquad
 p=\frac{1}{\sqrt{\st}}M\qt,\quad r=\frac{1}{\sqrt{s}}Mq.
\end{equation}
The vectors satisfy $p\geew|r|$, and consequently $\Delta\geew\zero$. If the Kelmans transformation is not proper, then $\Delta=\zero$. If it is proper, then $\Delta$ has rank two and inertia $(1,1,d-2)$, where the entries count positive, negative, and zero eigenvalues.
\end{theorem}
\end{samepage}
\begin{proof}
Apply~\eqref{eq:block-inverse} to the two matrices $xI-Q^TAQ$ and $xI-Q^T\widetilde AQ$. Their leading inverse blocks are $M+s^{-1}Mqq^TM$ and $M+\st^{-1}M\qt\qt^TM$, respectively. Subtraction gives~\eqref{eq:rank-two}. Since $M\geew\zero$, we have $M\qt\geew|Mq|$; together with $0<\st\leq s$, this gives $p\geew|r|$. Hence $p_ip_j-r_ir_j\geq0$ for all $i,j$. If $q$ has one sign, then $\qt=\pm q$, the Schur complements agree, and $\Delta=\zero$.

If $q$ has both signs, $q$ and $\qt$ are linearly independent. Invertibility of $M$ implies independence of $p,r$. For a real vector $z$, the quadratic form is $(p^Tz)^2-(r^Tz)^2$. The orthogonal projection of $p$ onto $r^\perp$ gives a vector on which the form is positive; the projection of $r$ onto $p^\perp$ gives one on which it is negative. Both projections are nonzero by independence. Since $\rank\Delta\leq2$, the inertia is $(1,1,d-2)$.
\end{proof}

For a proper Kelmans transformation, the negative eigenvalue of $\Delta(x)$ makes its restriction to the coupling space relevant to the comparison. Theorem~\ref{thm:universal} gives the precise condition.

\subsection{Support and strictness}
Let $\mathcal I$ be the component of the support graph of $F$ containing index $1$, and let $T$ be the union of the components of the support graph of $A$ containing $u$ or $v$. Under the basis given by $P$, the indices in $T$ correspond to $\mathcal I$, with $u,v$ replaced by the single basis vector $e_+$.

\begin{samepage}
\begin{theorem}\label{thm:local-strictness}
Suppose the Kelmans transformation is proper. For every $x>\rho(\widetilde A)$,
\[
 \Delta_{ij}(x)>0\quad\Longleftrightarrow\quad i,j\in\mathcal I;
 \qquad \Delta_{ij}(x)=0\quad\text{otherwise}.
\]
Moreover,
\[
 \rho(\widetilde A[T])>\rho(A[T]).
\]
If $c=\rho(A[T^c])$, with $c=0$ when $T^c$ is empty, then
\[
 \rho(\widetilde A)=\rho(A)\quad\Longleftrightarrow\quad c\geq\rho(\widetilde A[T]).
\]
\end{theorem}
\end{samepage}
\begin{proof}
For $2\leq j\leq d$, $F_{1j}\geq|q_j|$. Thus the support of $q$ is contained in $\mathcal I$. Lemma~\ref{lem:nonnegative-resolvent} shows that $M$ is block diagonal over the components of $F$ and that $M[\mathcal I]$ is strictly positive entrywise. Formula~\eqref{eq:rank-two} therefore gives $p_i=r_i=0$ for $i\notin\mathcal I$. For $i\in\mathcal I$, every coefficient $M_{ij}$ with $j\in\mathcal I$ is positive. Since $q$ has both signs, the strict triangle inequality gives $(M\qt)_i>|(Mq)_i|$, and hence $p_i>|r_i|$. Thus $p_ip_j>|r_i r_j|\geq r_i r_j$ for $i,j\in\mathcal I$, proving the support formula.

For the local spectral-radius inequality, restrict~\eqref{eq:normal-form} to the basis vectors $e_+$, $e_w$ with $w\in T\setminus\{u,v\}$, and $e_-$, and denote the resulting bordered matrices by $B_T,\widetilde B_T$; they are orthogonally similar to $A[T],\widetilde A[T]$, respectively. Choose $t_0\geq0$ so that both shifted matrices have nonnegative diagonal entries. Then $|B_T+t_0I|\leew\widetilde B_T+t_0I$, and $\widetilde B_T+t_0I$ is irreducible because $F[\mathcal I]$ is irreducible and $q\ne\zero$. Subtracting the common scalar $t_0$ from the Rayleigh-quotient inequality gives the following comparison. For a real unit eigenvector $y=\colvec zt$ associated with $\lambda_{\max}(B_T)$,
\[
 \lambda_{\max}(B_T)=y^TB_Ty
 \leq |y|^T\widetilde B_T|y|\leq\lambda_{\max}(\widetilde B_T).
\]
If equality holds, the Perron--Frobenius theorem for the shifted $\widetilde B_T$ gives $|y_i|>0$ for every entry of $y$. Equality in the first inequality requires the signs of $z$ to agree along every positive edge of $F[\mathcal I]$, hence throughout $\mathcal I$. It also requires $q_i z_i t>0$ for every nonzero $q_i$, which is impossible because $q$ has both signs. This proves strictness. Finally, the blocks on $T^c$ are unchanged and disconnected from $T$. Taking the maximum of the two component spectral radii gives the global equality criterion.
\end{proof}

Nonnegative vectors with equal entries at $u,v$ belong to $\Ran P$ and have nonnegative coefficients in the basis given by $P$. Theorems~\ref{thm:increment} and~\ref{thm:local-strictness} therefore give the following bilinear resolvent comparison and its strictness criterion.

\begin{samepage}
\begin{corollary}\label{cor:weighted-resolvent}
Let $a,c\in\R^n$ be nonnegative, with $a_u=a_v$ and $c_u=c_v$. Then
\[
 a^T(xI-\widetilde A)^{-1}c\geq a^T(xI-A)^{-1}c
 \qquad(x>\rho(\widetilde A)).
\]
The inequality is strict if and only if the Kelmans transformation is proper and each of $a,c$ has a positive entry indexed by $T$.
\end{corollary}
\end{samepage}
\begin{proof}
Let $\widehat a=P^Ta$ and $\widehat c=P^Tc$, so that $a=P\widehat a$ and $c=P\widehat c$; the vectors $\widehat a,\widehat c$ are nonnegative. The difference is $\widehat a^T\Delta\widehat c$. Theorems~\ref{thm:increment} and~\ref{thm:local-strictness} give $\widehat a^T\Delta\widehat c\geq0$.
\end{proof}

For a graph $G$, its coronal is
\[
 \Gamma_G(x)=\one^T(xI-A(G))^{-1}\one.
\]
Corollary~\ref{cor:weighted-resolvent} shows that a proper Kelmans transformation strictly increases $\Gamma_G(x)$ for every $x>\rho(\widetilde A)$, even if an unaffected component prevents strict increase of the spectral radius.

\begin{example}\label{ex:two-edges}\label{rem:cone-needed}
Let $G=2K_2$, with edges $u\ell$ and $vk$. The graph Kelmans transformation toward $u$ from $v$ gives $K_{1,2}$ and an isolated vertex. Thus the local spectral radius increases from $1$ to $\sqrt2$. With $k<\ell$, the columns of $P$ are $(e_+,e_k,e_\ell)$, and
\[
 F=\begin{bmatrix}0&1/\sqrt2&1/\sqrt2\\1/\sqrt2&0&0\\1/\sqrt2&0&0\end{bmatrix},
 \qquad q=[0,-1/\sqrt2,1/\sqrt2]^T,\qquad\delta=0.
\]
Here $F$ is irreducible, so every entry of $\Delta$ is positive; its inertia is $(1,1,1)$. Adding a direct summand $[5]$ leaves both global spectral radii equal to $5$. The local inequality is strict, and the unchanged summand determines the global spectral radius.

The sign and equal-coordinate conditions in Corollary~\ref{cor:weighted-resolvent} can be seen in two examples. If $G$ consists of the edge $vw$ and the isolated vertex $u$, then the unequal-coordinate vector $a=e_v$ gives
\[
 a^T\bigl((xI-\widetilde A)^{-1}-(xI-A)^{-1}\bigr)a=-\frac1{x(x^2-1)}<0\qquad(x>1).
\]
For the original graph $2K_2$, the signed vector $a=e_k-e_\ell$ has $a_u=a_v=0$, and the corresponding difference equals $-2/[x(x^2-1)]$ for $x>\sqrt2$.
\end{example}
\section{Comparison of Hermitian block extensions}\label{sec:block-comparison}

Let $H_j\in\C^{n_j\times n_j}$ be Hermitian matrices and let $B_j\in\C^{n_j\times m}$, for $j=0,1$, with $m\geq1$. The orders $n_0,n_1$ may be different. Define
\[
 \begin{gathered}
 E_j(C)=\begin{bmatrix}H_j&B_j\\B_j^*&C\end{bmatrix},\qquad h_j=\lambda_{\max}(H_j),\\
 \mathcal R_j(x)=B_j^*(xI-H_j)^{-1}B_j\quad(x>h_j),
 \end{gathered}
\]
where $C$ is an arbitrary Hermitian matrix of order $m$. We compare the extensions using $h_0,h_1$ and the compressed resolvents $\mathcal R_0,\mathcal R_1$.

\begin{samepage}
\begin{theorem}\label{thm:universal}
The following conditions are equivalent.
\begin{enumerate}
\raggedright
\item For every Hermitian matrix $C$, $\lambda_{\max}(E_1(C))\geq\lambda_{\max}(E_0(C))$.
\item $h_1\geq h_0$, and $\mathcal R_1(x)\succeq\mathcal R_0(x)$ for every $x>h_1$.
\item $h_1\geq h_0$, and for every $x>h_1$,
\[
 \{C:\lambda_{\max}(E_1(C))\leq x\}\subseteq\{C:\lambda_{\max}(E_0(C))\leq x\}.
\]
\item $h_1\geq h_0$, and $n_+(E_1(C)-xI)\geq n_+(E_0(C)-xI)$ for every Hermitian matrix $C$ and every $x>h_1$.
\end{enumerate}
\end{theorem}
\end{samepage}
\begin{proof}
For $x>\max\{h_0,h_1\}$, both $xI-H_j$ are positive definite. Applying Lemma~\ref{lem:schur} to $xI-E_j(C)$ gives
\begin{align*}
 \lambda_{\max}(E_j(C))\leq x
 &\quad\Longleftrightarrow\quad C\preceq xI-\mathcal R_j(x),\\
 n_+(E_j(C)-xI)&=n_-(xI-C-\mathcal R_j(x)).
\end{align*}
If (2) holds, then $xI-\mathcal R_1(x)\preceq xI-\mathcal R_0(x)$, so the first equivalence implies that condition~(3) holds. Also $xI-C-\mathcal R_1(x)\preceq xI-C-\mathcal R_0(x)$; eigenvalue monotonicity~\eqref{eq:variational-facts} and the displayed inertia identity imply condition~(4). Conversely, under (3), the choice $C=xI-\mathcal R_1(x)$ belongs to the first set in (3). Condition~(3) then implies $C\preceq xI-\mathcal R_0(x)$, and hence $\mathcal R_1(x)\succeq\mathcal R_0(x)$.

If (1) fails, choose $\lambda_{\max}(E_1(C))<x<\lambda_{\max}(E_0(C))$. The interlacing lemma~\eqref{eq:interlacing} gives $x>h_1$. This $C$ lies in the first set of (3), but not the second, contradicting (3). At the same value of $x$, $E_1(C)-xI$ has no positive eigenvalue and $E_0(C)-xI$ has at least one. Thus (4) also implies (1).

It remains to derive (2) from (1). Taking $C=-tI$ and letting $t\to\infty$ gives $\lambda_{\max}(E_j(-tI))\to h_j$. Indeed, the interlacing lemma gives the lower bound $h_j$, and the largest eigenvalue of
\[
 \begin{bmatrix}h_j&\|B_j\|\\\|B_j\|&-t\end{bmatrix}
\]
gives an upper bound, because for a unit vector $\colvec yz$ the quadratic form of $E_j(-tI)$ is at most
\[
 h_j\|y\|^2+2\|B_j\|\|y\|\|z\|-t\|z\|^2.
\]
The largest eigenvalue of this $2\times2$ matrix is
\[
 \frac{h_j-t+\sqrt{(h_j+t)^2+4\|B_j\|^2}}2\longrightarrow h_j.
\] Thus $h_1\geq h_0$. Condition (1) also gives the set inclusion in (3): if $\lambda_{\max}(E_1(C))\leq x$, then $\lambda_{\max}(E_0(C))\leq x$. The preceding argument now gives the resolvent order in (2).
\end{proof}

\phantomsection\label{ex:internal-order}
The internal condition compares the two largest eigenvalues. For example, take $H_0=\diag(1,0,-2)$, $H_1=\diag(2,-10,-3)$, and $B_0=B_1=e_1b^*$ for a nonzero $b\in\C^m$. Then $h_1=2>1=h_0$, and
\[
 \mathcal R_1(x)-\mathcal R_0(x)=\frac{1}{(x-2)(x-1)}bb^*\succeq\zero\qquad(x>2).
\]
Theorem~\ref{thm:universal} gives comparison for every $C$ in this example, where $H_1-H_0$ is indefinite.

Condition~(4) also compares the ordered eigenvalues above $h_1$:
\[
 \lambda_k(E_0(C))>h_1
 \quad\Longrightarrow\quad
 \lambda_k(E_1(C))\geq\lambda_k(E_0(C)).
\]
Indeed, let $x$ increase to $\lambda_k(E_0(C))$ from below in~(4). The interlacing lemma gives $k\leq m$, so the indicated eigenvalue of $E_1(C)$ exists even when the internal orders are different.

A negative eigenvalue of the compressed resolvent difference gives common blocks for which the eigenvalue inequality is reversed. At a fixed spectral parameter, the next proposition determines an open ball of such blocks and its largest possible radius.

\begin{samepage}
\begin{proposition}\label{prop:reversing-ball}
Fix $x>\max\{h_0,h_1\}$ and suppose $\lambda_{\min}(\mathcal R_1(x)-\mathcal R_0(x))=-\eta<0$.
Let $C_*=xI-\mathcal R_1(x)-\frac{\eta}{2}I$. For every Hermitian matrix $K$ with $\|K\|<\eta/2$,
\[
 \lambda_{\max}(E_0(C_*+K))>x>\lambda_{\max}(E_1(C_*+K)).
\]
Moreover, $\eta/2$ is the largest radius of an open operator-norm ball of Hermitian common blocks for which these two strict inequalities hold throughout the ball.
\end{proposition}
\end{samepage}
\begin{proof}
Let $D=\mathcal R_1(x)-\mathcal R_0(x)$. The two Schur complements of $xI-E_j(C_*+K)$ are
\[
 S_1=\frac\eta2I-K\succ\zero,\qquad S_0=D+\frac\eta2I-K.
\]
For a unit vector $y$ with $Dy=-\eta y$, one has $y^*S_0y\leq-\eta/2+\|K\|<0$. Thus $xI-E_1$ is positive definite, whereas $xI-E_0$ has a negative eigenvalue.

For optimality, let $L_j=xI-\mathcal R_j(x)$ and suppose an open ball with center $C_0$ and radius $r>0$ satisfies the two strict inequalities of the proposition. Schur complementation gives $C\prec L_1$ and $\lambda_{\max}(C-L_0)>0$ throughout that ball. Substituting $C=C_0+tI$ and $C=C_0-tI$ for $0<t<r$, and then letting $t\uparrow r$, gives
\[
 C_0\preceq L_1-rI,\qquad \lambda_{\max}(C_0-L_0)\geq r.
\]
Consequently,
\[
 r\leq\lambda_{\max}(C_0-L_0)
 \leq\lambda_{\max}(L_1-L_0)-r=\eta-r,
\]
which proves $r\leq\eta/2$.
\end{proof}

For a coupling with zero border row, subtraction of the two compressed resolvents cancels the common term and leaves a difference of rank-one matrices.

\begin{samepage}
\begin{corollary}\label{cor:common-border}
Let
\[
 H_j=\begin{bmatrix}F&q_j\\q_j^*&\delta_j\end{bmatrix},\qquad B_j=\colvec Z{\zero},
\]
where $F$ is a Hermitian matrix, $\delta_j\in\R$, and $h_1\geq h_0$. For $x>h_1$, let
\[
 M=(xI-F)^{-1},\qquad s_j=x-\delta_j-q_j^*Mq_j,\qquad
 v_j=\frac{1}{\sqrt{s_j}}Z^*Mq_j.
\]
Then the comparison in Theorem~\ref{thm:universal}(1) holds if and only if, for every $x>h_1$, $v_0(x)=\theta(x)v_1(x)$ for some $|\theta(x)|\leq1$.
\end{corollary}
\end{samepage}
\begin{proof}
The interlacing lemma gives $xI-F\succ\zero$, and Lemma~\ref{lem:schur} gives $s_j>0$ for $j=0,1$. The inverse formula~\eqref{eq:block-inverse} yields
\[
 \mathcal R_j(x)=Z^*\left(M+\frac1{s_j}Mq_jq_j^*M\right)Z
 =Z^*MZ+v_jv_j^*.
\]
Subtracting cancels $Z^*MZ$. Lemma~\ref{lem:rank-one-order} characterizes positivity of the remaining difference $v_1v_1^*-v_0v_0^*$, which is the resolvent condition in Theorem~\ref{thm:universal}.
\end{proof}
\section{Nonnegative rank-one couplings}\label{sec:rank-one}

We now return to $A,\widetilde A$ related by a Kelmans transformation. Let $a\in\R^n$ be nonnegative with $a_u=a_v$, let $b\in\C^m$, and let $C\in\C^{m\times m}$ be a Hermitian matrix. For a real parameter $\tau$, define
\[
 L_\tau(A)=\begin{bmatrix}\tau A&ab^*\\ba^T&C\end{bmatrix}.
\]
We keep $a,b,C$ fixed when $A$ is replaced by $\widetilde A$.

\begin{samepage}
\begin{theorem}\label{thm:rank-one-monotonicity}
The inequality
\begin{equation}\label{eq:rank-one-monotonicity}
 \lambda_{\max}(L_\tau(\widetilde A))\geq\lambda_{\max}(L_\tau(A))
\end{equation}
holds for every allowed $A,u,v,a,b,C$ if and only if $\tau\geq0$.
\end{theorem}
\end{samepage}
\begin{proof}
For $\tau=0$ the internal blocks are both zero, so $L_0(\widetilde A)=L_0(A)$. For $\tau>0$, Lemma~\ref{lem:inverse-comparison} gives $\lambda_{\max}(\tau \widetilde A)\geq\lambda_{\max}(\tau A)$. For $x>\tau\rho(\widetilde A)$, the compressed resolvent difference is
\[
 \frac1\tau\left[a^T\left(\left(\frac x\tau I-\widetilde A\right)^{-1}
 -\left(\frac x\tau I-A\right)^{-1}\right)a\right]bb^*\succeq\zero
\]
by Corollary~\ref{cor:weighted-resolvent}. Theorem~\ref{thm:universal} proves the comparison.

For necessity, let $\tau=-c<0$ and let $A$ be the adjacency matrix of the path with edges $12,23,14$. The graph Kelmans transformation toward $u=3$ from $v=4$ gives $K_3\cup K_1$. If $\varphi=(1+\sqrt5)/2$, then $\lambda_{\max}(-cA)=c\varphi$ and $\lambda_{\max}(-c\widetilde A)=c$.
To retain a nonzero positive coupling, take $a=\varepsilon\one_4$, $b=1$, and $C=0$. The off-diagonal perturbation has norm $2\varepsilon$. The eigenvalue bound in~\eqref{eq:variational-facts} gives
\[
 \lambda_{\max}(L_{-c}(A))\geq c\varphi-2\varepsilon
 >c+2\varepsilon\geq\lambda_{\max}(L_{-c}(\widetilde A))
\]
whenever $0<\varepsilon<c(\varphi-1)/4$.
\end{proof}

\subsection{Equality at and above the internal largest eigenvalue}
Equality in Theorem~\ref{thm:rank-one-monotonicity} can be caused by an eigenvector supported on the complementary block. If $Cy=\mu y$ and $b^*y=0$ for a nonzero $y$, then $\colvec{\zero}y$ is an eigenvector with eigenvalue $\mu$ before and after the transformation. This eigenvalue may be largest in both matrices. The generalized Schur complement gives an equality criterion valid also at $\mu=\tau\rho(\widetilde A)$, where the leading block is singular. Throughout this subsection, assume $a\geew\zero$ and $a_u=a_v$.

\begin{samepage}
\begin{proposition}\label{prop:equality}
Suppose $\tau>0$ and let $\mu=\lambda_{\max}(L_\tau(A))$. Equality in~\eqref{eq:rank-one-monotonicity} holds if and only if, with $K=\mu I-\tau \widetilde A$,
\begin{equation}\label{eq:equality-general}
 K\succeq\zero,\qquad (I-KK^\dagger)ab^*=\zero,\qquad
 \mu I-C-(a^TK^\dagger a)bb^*\succeq\zero.
\end{equation}
If $\mu>\tau\rho(\widetilde A)$, define
\[
 \gamma=a^T(\mu I-\tau A)^{-1}a,\quad
 \widetilde\gamma=a^T(\mu I-\tau \widetilde A)^{-1}a,\quad
 S=\mu I-C-\gamma bb^*,\quad\eta=\widetilde\gamma-\gamma.
\]
Then $S\succeq\zero$ is singular and $\eta\geq0$. Equality holds if $\eta=0$; if $\eta>0$, it holds exactly when $b\in\Ran S$ and $\eta b^*S^\dagger b\leq1$.
\end{proposition}
\end{samepage}
\begin{proof}
Theorem~\ref{thm:rank-one-monotonicity} already gives $\lambda_{\max}(L_\tau(\widetilde A))\geq\mu$. Equality is therefore equivalent to $\mu I-L_\tau(\widetilde A)\succeq\zero$. Formula~\eqref{eq:generalized-schur}, with $X=ab^*$ and $D=\mu I-C$, gives~\eqref{eq:equality-general}. If $\mu>\tau\rho(\widetilde A)$, both leading blocks are positive definite, with Schur complements $S$ and $S-\eta bb^*$. Corollary~\ref{cor:weighted-resolvent} at $x=\mu/\tau$ gives $\eta\geq0$. Since $\mu I-L_\tau(A)$ is positive semidefinite and singular, Lemma~\ref{lem:schur} shows that $S$ is positive semidefinite and singular. When $\eta=0$ the two Schur complements agree. Otherwise, the rank-one criterion~\eqref{eq:downdate} gives the stated condition.
\end{proof}

The preceding equality criterion, together with strictness of the weighted resolvent comparison, gives a sufficient condition for strict increase of the largest eigenvalue. For the following corollary, retain $\mu=\lambda_{\max}(L_\tau(A))$ from Proposition~\ref{prop:equality}.

\begin{samepage}
\begin{corollary}\label{cor:rank-one-strictness}
Suppose $A$ is irreducible, the Kelmans transformation is proper, $a\ne\zero$, $b\ne\zero$, and $\tau>0$. If $\lambda_{\max}(L_\tau(A))>\lambda_{\max}(C)$, then the inequality in~\eqref{eq:rank-one-monotonicity} is strict. In particular, equality can hold only if
\[
 \mu=\lambda_{\max}(C),\qquad
 \ker(\mu I-C)\cap b^\perp\ne\{\zero\}.
\]
\end{corollary}
\end{samepage}
\begin{proof}
If $\mu<\tau\rho(\widetilde A)$, the interlacing lemma gives strictness. If $\mu=\tau\rho(\widetilde A)$, consider the range condition in~\eqref{eq:equality-general}. After the transformation, every neighbor of $v$ other than $u$ is a neighbor of $u$. Replacing $v$ by $u$ in paths of the connected support graph of $A$ and deleting consecutive repetitions therefore puts every vertex except possibly $v$ in the component of $u$. An isolated $v$ has diagonal entry $\widetilde a_{vv}\leq \widetilde a_{uu}$, while the component of $u$ has a positive edge. For a positive Perron vector $z$ of this component, each row has a positive off-diagonal contribution, so $\rho z_i=\widetilde a_{ii}z_i+\sum_{j\ne i}\widetilde a_{ij}z_j>\widetilde a_{ii}z_i$. In particular, the component has spectral radius greater than $\widetilde a_{uu}$. Thus a Perron eigenvector of $\widetilde A$ is positive in every entry except possibly the one indexed by $v$. The conditions $a\geew\zero$, $a\ne\zero$, and $a_u=a_v$ imply a positive inner product with this eigenvector. If $z$ denotes this Perron vector extended by zero, then $z^*a>0$ and $z^*ab^*\ne\zero$. Proposition~\ref{prop:equality} therefore gives
$\lambda_{\max}(L_\tau(\widetilde A))>\mu$.

It remains to consider $\mu>\tau\rho(\widetilde A)$. Corollary~\ref{cor:weighted-resolvent} gives $\eta>0$. If equality holds, Proposition~\ref{prop:equality} gives $b\perp\ker S$. For any nonzero $y\in\ker S$, one therefore has $(\mu I-C)y=\zero$ and $b^*y=0$. The interlacing lemma gives $\mu\geq\lambda_{\max}(C)$, while this eigenvector gives the reverse inequality. This proves the two necessary conditions for equality and the strictness assertion.
\end{proof}

For example, take $C=\diag(L,C_0)$ and $b=\colvec0{b_0}$. If $L$ exceeds the largest eigenvalues of both remaining coupled blocks, then both full matrices have largest eigenvalue $L$. In this case, the unchanged eigenvector supported on the first coordinate of $C$ gives equality for the full matrices.
\section{Coupling subspaces and polynomial criteria}\label{sec:subspaces}

For a nonnegative rank-one coupling, entrywise resolvent comparison gives the required scalar inequality. For a complex coupling matrix, the resolvent difference must be positive semidefinite on its column space. Retain $P,F,q,\qt$ from~\eqref{eq:normal-form} and $M,s,\st,\Delta,p,r$ from~\eqref{eq:schur-data}--\eqref{eq:rank-two}. In this section, subspaces are complex and orthogonality refers to the standard Hermitian inner product. For $Z\in\C^{d\times m}$ and $\tau>0$, define
\[
 N_\tau(A;Z,C)=\begin{bmatrix}\tau A&PZ\\Z^*P^T&C\end{bmatrix}.
\]
We call $\U=\Ran Z$ an \emph{admissible coupling subspace} if
\[
 \lambda_{\max}(N_\tau(\widetilde A;Z,C))\geq\lambda_{\max}(N_\tau(A;Z,C))
\]
for every Hermitian matrix $C$ of order $m$. The following criterion expresses admissibility in terms of $\U$ and applies for every $\tau>0$.

\begin{samepage}
\begin{theorem}\label{thm:coupling-criterion}
The subspace $\U=\Ran Z$ is admissible if and only if
\begin{equation}\label{eq:compressed-order}
 Z^*\Delta(x)Z\succeq\zero\qquad\text{for all }x>\rho(\widetilde A).
\end{equation}
Equivalently, for each such $x$ there is $\theta(x)\in\C$, $|\theta(x)|\leq1$, such that
\begin{equation}\label{eq:pointwise-collinearity}
 Z^*r(x)=\theta(x)Z^*p(x).
\end{equation}
For a proper Kelmans transformation, every admissible subspace has dimension at most $d-1$.
\end{theorem}
\end{samepage}
\begin{proof}
Lemma~\ref{lem:inverse-comparison} gives $\tau\rho(\widetilde A)\geq\tau\rho(A)$. The resolvent condition in Theorem~\ref{thm:universal} is therefore $\tau^{-1}Z^*\Delta(x/\tau)Z\succeq\zero$ for all $x>\tau\rho(\widetilde A)$. Substituting $x=\tau t$ gives~\eqref{eq:compressed-order}. By~\eqref{eq:rank-two}, the compressed matrix is $(Z^*p)(Z^*p)^*-(Z^*r)(Z^*r)^*$, so Lemma~\ref{lem:rank-one-order} gives~\eqref{eq:pointwise-collinearity}. Finally, $y^*Z^*\Delta Zy=(Zy)^*\Delta(Zy)$ shows that the condition depends only on $\Ran Z$. For a proper Kelmans transformation, Theorem~\ref{thm:increment} gives a vector $y$ with $y^*\Delta y<0$. Hence every admissible subspace is a proper subspace of $\C^d$ and has dimension at most $d-1$.
\end{proof}

At each $x$, one scalar $\theta(x)$ gives~\eqref{eq:pointwise-collinearity} for every column of $Z$; equivalently, $\U\subseteq(r(x)-\theta(x)p(x))^\perp$. The factor can depend on $x$, as in Example~\ref{ex:nonconstant}. If the compressed difference has a negative eigenvalue at some $x$, Proposition~\ref{prop:reversing-ball} gives an open ball of common complementary matrices $C$ for which the largest-eigenvalue inequality is reversed.

For $m=1$, criterion~\eqref{eq:pointwise-collinearity} reduces to $|Z^*r|\leq|Z^*p|$ for every $x>\rho(\widetilde A)$, including $Z^*r=0$ whenever $Z^*p=0$. Generally, the condition compares the two projections on $\U$ through this proportionality and magnitude bound.

\begin{remark}
An admissible subspace may contain vectors with a nonzero projection onto the negative eigenspace of $\Delta(x)$. Admissibility requires the full quadratic form on the subspace to be nonnegative, as expressed by~\eqref{eq:pointwise-collinearity}.
\end{remark}

\subsection{Admissible hyperplanes}
For a proper Kelmans transformation, admissible subspaces have dimension at most $d-1$. We first determine when this bound can be attained. Let
\[
 \mathcal S=\spanop\{q,\qt\},\qquad V=[q\ \qt],
 \qquad \Pi_{\mathcal S}=V(V^*V)^{-1}V^*.
\]
Properness makes $q,\qt$ linearly independent, so $\Pi_{\mathcal S}$ is the orthogonal projector onto the two-dimensional subspace $\mathcal S$.

\begin{samepage}
\begin{theorem}\label{thm:hyperplanes}
Suppose the Kelmans transformation is proper. A hyperplane $\U\subset\C^d$ is admissible if and only if there are constants $\kappa,\xi_0,\xi_1\in\C$, with $|\kappa|\leq1$, such that
\begin{align}
 w&=q-\kappa\qt,\qquad \U=w^\perp,\qquad
 Fw=\xi_0q+\xi_1\qt,\label{eq:hyperplane-normal}\\
 s(x)|x-\xi_0|^2&-\st(x)|\kappa x+\xi_1|^2\geq0
 \qquad(x>\rho(\widetilde A)).\label{eq:hyperplane-sign}
\end{align}
Here $M,s,\st$ are the resolvent and Schur complements of~\eqref{eq:schur-data}. For such a hyperplane, any matrix $Z$ with range $\U$ satisfies
\begin{equation}\label{eq:hyperplane-factor}
 Z^*M(x)q=\frac{\kappa x+\xi_1}{x-\xi_0}Z^*M(x)\qt,
 \qquad x>\rho(\widetilde A),
\end{equation}
and the denominator is nonzero throughout this interval.
\end{theorem}
\end{samepage}
\begin{proof}
Choose $Z$ with orthonormal columns spanning $\U$. From $M=x^{-1}I+O(x^{-2})$ and $s,\st=x+O(1)$, the resolvent expansion at infinity gives
\[
 \lim_{x\to\infty}x^3\Delta(x)=\qt\qt^T-qq^T.
\]
If $\U$ is admissible, its compression is positive semidefinite. Lemma~\ref{lem:rank-one-order} gives $Z^*q=\kappa Z^*\qt$ for some $|\kappa|\leq1$. Since $q,\qt$ are independent, $w=q-\kappa\qt$ is nonzero and spans $\U^\perp$.

For each $x>\rho(\widetilde A)$, the vectors $Mq,M\qt$ are independent, but their projections onto $\U$ are dependent. A nonzero linear combination therefore lies in the one-dimensional kernel $\U^\perp$ of $Z^*$, giving $w\in\spanop\{Mq,M\qt\}$. Multiplication by $xI-F$ gives $xw-Fw\in\mathcal S$. Since $w\in\mathcal S$, also $Fw\in\mathcal S$, with unique coefficients $\xi_0,\xi_1$ as in~\eqref{eq:hyperplane-normal}. Applying $Z^*M$ to
\[
 (xI-F)w=(x-\xi_0)q-(\kappa x+\xi_1)\qt
\]
gives $(x-\xi_0)Z^*Mq=(\kappa x+\xi_1)Z^*M\qt$. In fact, $Z^*M\qt$ is nonzero: otherwise admissibility implies $Z^*Mq=\zero$, putting two independent vectors in the one-dimensional space $\U^\perp$. If $x-\xi_0=0$, the displayed relation consequently gives $\kappa x+\xi_1=0$, contradicting invertibility of $xI-F$. Dividing by $x-\xi_0$ gives~\eqref{eq:hyperplane-factor}.
Let $v=Z^*M\qt$. Substituting~\eqref{eq:hyperplane-factor}
into~\eqref{eq:rank-two} yields
\[
 Z^*\Delta Z=
 \left(\frac{1}{\st}
 -\frac{|\kappa x+\xi_1|^2}{s|x-\xi_0|^2}\right)vv^*.
\]
Since $\U$ is admissible, $Z^*\Delta Z\succeq\zero$.
The vector $v$ is nonzero, so the scalar coefficient is nonnegative.
Multiplying by $s\st|x-\xi_0|^2>0$
gives~\eqref{eq:hyperplane-sign}.

Conversely, suppose~\eqref{eq:hyperplane-normal}--\eqref{eq:hyperplane-sign} hold. A zero of $x-\xi_0$ in the interval implies $\kappa x+\xi_1=0$, because $\st>0$. This would give $(xI-F)w=\zero$, a contradiction. Relation~\eqref{eq:hyperplane-factor} therefore holds throughout the interval. For $v=Z^*M\qt$, it gives
\[
 Z^*\Delta Z=
 \frac{s|x-\xi_0|^2-\st|\kappa x+\xi_1|^2}
 {s\st|x-\xi_0|^2}\,vv^*\succeq\zero.
\]
Theorem~\ref{thm:coupling-criterion} proves admissibility.
\end{proof}

The condition $Fw\in\mathcal S$ gives the following dimension bound.

\begin{samepage}
\begin{corollary}\label{cor:dimension-obstruction}
For a proper Kelmans transformation with $\rank((I-\Pi_{\mathcal S})FV)=2$, every admissible coupling subspace has dimension at most $d-2$.
\end{corollary}
\end{samepage}
\begin{proof}
Let $w=V\alpha\in\mathcal S$. The condition $Fw\in\mathcal S$ is equivalent to $(I-\Pi_{\mathcal S})FV\alpha=\zero$. By the rank assumption, this gives $\alpha=\zero$, hence $w=\zero$. Theorem~\ref{thm:hyperplanes} therefore excludes dimension $d-1$, and Theorem~\ref{thm:coupling-criterion} gives dimension at most $d-2$.
\end{proof}

To test admissibility of a hyperplane, its normal vector $w=q-\kappa\qt$ must first satisfy $|\kappa|\leq1$ and $[1,-\kappa]^T\in\ker((I-\Pi_{\mathcal S})FV)$. For each such normal vector, determine $\xi_0,\xi_1$ from~\eqref{eq:hyperplane-normal} and check~\eqref{eq:hyperplane-sign}. The rational factor has numerator and denominator of degree at most one, regardless of $d$.

\subsection{Clearing denominators}
Define the real polynomials and polynomial matrix
\[
 \begin{aligned}
 \pi(x)&=\det(xI-F),\qquad J(x)=\adj(xI-F),\\
 \beta(x)&=(x-\delta)\pi(x)-q^TJ(x)q,\\
 \bt(x)&=(x-\delta)\pi(x)-\qt^TJ(x)\qt.
 \end{aligned}
\]
On $x>\rho(\widetilde A)$, the quantities $\pi$, $\beta=\pi s$, and $\bt=\pi\st$ are positive. Let
\[
 X(x)=Z^*J(x)\qt,\qquad Y(x)=Z^*J(x)q.
\]
Norm squares of complex polynomial vectors mean evaluation at real $x$ followed by the ordinary squared norm. They are real polynomials in $x$.

The proportionality requirement is now expressed by polynomial minors, and the magnitude bound becomes one scalar polynomial inequality.

\begin{samepage}
\begin{theorem}\label{thm:polynomial-test}
The subspace $\Ran Z$ is admissible if and only if
\begin{align}
 X_iY_j-X_jY_i&\equiv0\quad(1\leq i<j\leq m),\label{eq:minors}\\
 \beta(x)\|X(x)\|^2-\bt(x)\|Y(x)\|^2&\geq0\quad(x>\rho(\widetilde A)).\label{eq:sign-test}
\end{align}
The degrees in~\eqref{eq:minors} are at most $2d-2$, and the degree in~\eqref{eq:sign-test} is at most $3d-1$. If the entries of $A$ and $Z$ are algebraic, these conditions can be decided by finitely many exact polynomial operations and real-root isolation.
\end{theorem}
\end{samepage}
\begin{proof}
Equation~\eqref{eq:rank-two} gives
\[
 Z^*\Delta Z=\frac{1}{\pi\bt}XX^*-\frac{1}{\pi\beta}YY^*.
\]
By Lemma~\ref{lem:rank-one-order}, this is positive semidefinite precisely when $X,Y$ are collinear and~\eqref{eq:sign-test} holds. The assertion includes $X=\zero$, since the scalar inequality then implies $Y=\zero$. Pointwise collinearity on an open interval is equivalent to the polynomial identities~\eqref{eq:minors}. The degree bounds follow from $\deg J\leq d-1$ and $\deg\beta,\deg\bt\leq d+1$. For algebraic data, the endpoint $\rho(\widetilde A)$ is algebraic. Isolating the real roots of the sign polynomial partitions the interval into finitely many intervals on which its sign is constant; the identically zero case is immediate.
\end{proof}

For a fixed $Z$, the remaining scalar inequality also has a semidefinite representation. We state the standard univariate positivity result in the form used here; see~\cite{MarshallSOS} for background on sums of squares.

\begin{samepage}
\begin{proposition}\label{prop:sos-certificate}
Let $\sigma\in\R[x]$ have degree $\nu\geq1$, let $h\in\R$, and let $z_k(t)=[1,t,\ldots,t^k]^T$. Then $\sigma(x)\geq0$ for all $x>h$ if and only if there are real symmetric positive semidefinite matrices $G_0,G_1$ such that
\[
 \sigma(h+t)=z_{\lfloor\nu/2\rfloor}(t)^TG_0z_{\lfloor\nu/2\rfloor}(t)
 +t\,z_{\lfloor(\nu-1)/2\rfloor}(t)^TG_1z_{\lfloor(\nu-1)/2\rfloor}(t).
\]
For a constant polynomial, nonnegativity is represented by $G_0=[\sigma]$ with no $G_1$ term; this includes the zero polynomial. Thus, for fixed $A,Z$, the sign condition in~\eqref{eq:sign-test} is equivalent to finitely many linear coefficient equations and at most two positive semidefinite constraints.
\end{proposition}
\end{samepage}
\begin{proof}
Continuity extends nonnegativity to $x=h$. The polynomial $\sigma(h+y^2)$ is therefore nonnegative on the real line. Factoring into real linear and irreducible quadratic factors shows that it is a sum of squares $\sum_j r_j(y)^2$, with $\deg r_j\leq\nu$: real roots have even multiplicity, each positive quadratic is a sum of two squares, and products are handled by $(a^2+b^2)(c^2+d^2)=(ac-bd)^2+(ad+bc)^2$. Write $r_j(y)=a_j(y^2)+y b_j(y^2)$. Substitute these expressions into the sum-of-squares identity and average the resulting identities at $y$ and $-y$. It follows that
\[
 \sigma(h+t)=\sum_j a_j(t)^2+t\sum_j b_j(t)^2,
 \qquad \deg a_j\leq\lfloor\nu/2\rfloor,\quad
 \deg b_j\leq\lfloor(\nu-1)/2\rfloor.
\]
The sums of outer products of the coefficient vectors of $a_j,b_j$ give $G_0,G_1$. Conversely, factoring positive semidefinite Gram matrices gives sums of squares and hence nonnegativity for $t\geq0$. Coefficient matching is linear in their entries.
\end{proof}

For prescribed $A,Z$, this gives a finite semidefinite feasibility problem, also when $Z$ is complex. An exact certificate consists of Gram matrices satisfying the coefficient equations and positive semidefiniteness.

\subsection{A description by coefficient kernels}
For polynomials $f,g\in\C[x]$ with $g\ne0$, set
\[
 D_{f,g}(x)=g(x)J(x)q-f(x)J(x)\qt
 =\sum_{k=0}^{N}d_kx^k,\qquad
 \mathsf D_{f,g}=[d_0\ \cdots\ d_N],
\]
and define the real polynomial
\[
 \Psi_{f,g}(x)=\beta(x)|g(x)|^2-\bt(x)|f(x)|^2.
\]
The space $\ker\mathsf D_{f,g}^*$ consists of vectors orthogonal to every coefficient $d_k$; equivalently, they are orthogonal to $D_{f,g}(x)$ for every real $x$. The next theorem combines this coefficient condition with scalar nonnegativity to characterize all admissible subspaces.

\begin{samepage}
\begin{theorem}\label{thm:subspace-classification}
A subspace $\U\subseteq\C^d$ is admissible if and only if there are coprime polynomials $f,g\in\C[x]$, $g\ne0$, with $\deg f,\deg g\leq d-1$ (allowing $f=0$), such that
\[
 \U\subseteq\ker\mathsf D_{f,g}^*,\qquad
 \Psi_{f,g}(x)\geq0\quad(x>\rho(\widetilde A)).
\]
For each fixed pair satisfying the scalar inequality, the largest subspace satisfying its polynomial relation is $\ker\mathsf D_{f,g}^*$, of dimension $d-\rank\mathsf D_{f,g}$. If the Kelmans transformation is proper, this dimension is at most $d-1$.
\end{theorem}
\end{samepage}
\begin{proof}
Choose $Z$ with range $\U$. Suppose first that $\U$ is admissible. If $X\equiv\zero$, Theorem~\ref{thm:polynomial-test} implies $Y\equiv\zero$, and $f=0,g=1$ works. Otherwise, the identities~\eqref{eq:minors} imply that $Y=(f/g)X$ over $\C(x)$ for some reduced fraction $f/g$. Since $gY=fX$ and $f,g$ are coprime, $g$ divides every entry of $X$. Hence
\[
 X=gV,\qquad Y=fV
\]
for a nonzero polynomial vector $V$. A nonzero entry of $V$ gives the bounds $\deg f,\deg g\leq d-1$. The identity $gY-fX=\zero$ is equivalent to $\U\subseteq\ker\mathsf D_{f,g}^*$. The sign polynomial in~\eqref{eq:sign-test} equals $\|V\|^2\Psi_{f,g}$. A nonzero polynomial vector is zero at only finitely many real points, so continuity proves the required inequality for $\Psi_{f,g}$ on the entire interval.

Conversely, suppose the two conditions of the theorem hold. A real zero $x_0>\rho(\widetilde A)$ of $g$ would not be a zero of $f$, by coprimality, and would give $\Psi_{f,g}(x_0)=-\bt(x_0)|f(x_0)|^2<0$. Thus $g$ has no zero there. The coefficient condition gives $Y=(f/g)X$, and
\[
 \beta\|X\|^2-\bt\|Y\|^2
 =\frac{\|X\|^2}{|g|^2}\Psi_{f,g}\geq0.
\]
Theorem~\ref{thm:polynomial-test} proves admissibility.

The assertion for a fixed pair is the coefficient identity itself. If the Kelmans transformation is proper and $D_{f,g}\equiv\zero$, multiplication by $xI-F$ gives $g\pi q-f\pi\qt=\zero$, contrary to the linear independence of $q,\qt$. Thus the coefficient matrix has positive rank.
\end{proof}

Suppose $A$ has real algebraic entries and the prescribed coprime pair $f,g$, with $g\ne0$, has algebraic coefficients, all effectively represented. Decide the scalar sign condition by real-root isolation. If it holds, form $\mathsf D_{f,g}$ and compute a basis of $\ker\mathsf D_{f,g}^*$. An admissible coupling of dimension $k$ with this relation then exists exactly when $\rank\mathsf D_{f,g}\leq d-k$. The full classification in Theorem~\ref{thm:subspace-classification} is given by varying the prescribed pair.

For a fixed pair $f,g$, Proposition~\ref{prop:sos-certificate} also applies to $\Psi_{f,g}$. The resulting feasibility problem is linear in the Gram matrices; choosing $f,g$ at the same time introduces nonlinear coefficient constraints. When $X\not\equiv\zero$, the quotient $f/g$ in the proof is unique. It describes the unnormalized relation $Z^*Mq=(f/g)Z^*M\qt$. The factor in~\eqref{eq:pointwise-collinearity} is instead
\[
 \theta(x)=\frac{f(x)}{g(x)}\sqrt{\frac{\st(x)}{s(x)}}.
\]
Thus the magnitude condition is $|f(x)/g(x)|\sqrt{\st(x)/s(x)}\leq1$.

For $f=\theta$ and $g=1$ with $|\theta|\leq1$, the scalar inequality follows from $\st\leq s$. The coefficient condition then has a cyclic-subspace form.

\begin{samepage}
\begin{corollary}\label{cor:krylov}
Fix $\theta\in\C$ with $|\theta|\leq1$. Every subspace orthogonal to
\begin{equation}\label{eq:krylov-space}
 \spanop\{F^k(q-\theta\qt):0\leq k<d\}
\end{equation}
is admissible.
\end{corollary}
\end{samepage}
\begin{proof}
Lemma~\ref{lem:resolvent-span} puts $M(q-\theta\qt)$ in~\eqref{eq:krylov-space}, so $Z^*Mq=\theta Z^*M\qt$. The normalized factor has modulus $|\theta|\sqrt{\st/s}\leq1$, proving admissibility by Theorem~\ref{thm:coupling-criterion}. The same lemma identifies~\eqref{eq:krylov-space} with the span of the coefficients of $J(x)(q-\theta\qt)$. Thus this construction is the case $f=\theta,g=1$ of Theorem~\ref{thm:subspace-classification}.
\end{proof}

\subsection{Reduction to the relevant eigenvalues of \texorpdfstring{$F$}{F}}
The spectral decomposition of $F$ gives a smaller common denominator for $Mq$ and $M\qt$. The degree bounds can therefore be expressed using the eigenvalues with nonzero projections of $q$ or $\qt$.

\begin{samepage}
\begin{proposition}\label{prop:reduced-degrees}
Suppose $q\ne\zero$ and write
$F=\sum_\lambda\lambda E_\lambda$,
where $\lambda$ ranges over the distinct eigenvalues of $F$
and $E_\lambda$ is the orthogonal projector onto its
$\lambda$-eigenspace. Let
\[
 \Lambda=\{\lambda:E_\lambda q\ne\zero
 \text{ or }E_\lambda\qt\ne\zero\},
 \qquad \ell=|\Lambda|,
 \qquad
 \chi(x)=\prod_{\lambda\in\Lambda}(x-\lambda).
\]
Define
\begin{align*}
 \Phi_0(x)
 &=\sum_{\lambda\in\Lambda}
   \frac{\chi(x)}{x-\lambda}E_\lambda q,
&
 \Phi_1(x)
 &=\sum_{\lambda\in\Lambda}
   \frac{\chi(x)}{x-\lambda}E_\lambda\qt,\\
 b_0(x)
 &=(x-\delta)\chi(x)-q^T\Phi_0(x),
&
 b_1(x)
 &=(x-\delta)\chi(x)-\qt^T\Phi_1(x).
\end{align*}
Then $\chi$ is the monic least common denominator
of the entries of $Mq$ and $M\qt$.
\end{proposition}
\end{samepage}

\begin{proof}
The spectral decomposition gives
\[
 Mq=\sum_{\lambda\in\Lambda}
       \frac{E_\lambda q}{x-\lambda},
 \qquad
 M\qt=\sum_{\lambda\in\Lambda}
       \frac{E_\lambda\qt}{x-\lambda}.
\]
For each $\lambda\in\Lambda$, at least one of the two
residues is nonzero. Thus each factor $x-\lambda$
occurs in the least common denominator, and all poles
are simple. This proves the assertion about $\chi$.

Multiplication by $\chi$ gives
\[
 \Phi_0=\chi Mq,\qquad \Phi_1=\chi M\qt,
\]
so both polynomial vectors have degree at most $\ell-1$.
Their definitions also give
\[
 b_0=\chi s,\qquad b_1=\chi\st,
 \qquad \deg b_0,\deg b_1\leq\ell+1.
\]
For $x>\rho(\widetilde A)$, we have $\chi>0$ and
$s,\st>0$, hence $b_0,b_1>0$.
\end{proof}

\begin{remark}\label{rem:reduced-criteria}
The polynomial test in Theorem~\ref{thm:polynomial-test}
takes the same form with
$Z^*\Phi_1,Z^*\Phi_0,b_0,b_1$
in place of $X,Y,\beta,\bt$, respectively.
The minors have degree at most $2\ell-2$,
and the scalar sign polynomial has degree at most $3\ell-1$.
In Theorem~\ref{thm:subspace-classification},
the coefficient polynomial and the scalar sign polynomial
can be written as
\[
 g\Phi_0-f\Phi_1,
 \qquad
 b_0|g|^2-b_1|f|^2,
\]
respectively, and coprime polynomials $f,g$ with
$\deg f,\deg g\leq\ell-1$ suffice.

To verify these statements, let $h=\pi/\chi$.
Since $\chi$ divides $\pi$, $h$ is a polynomial, and
$h(x)>0$ for $x>\rho(\widetilde A)$. Moreover,
\[
 Jq=h\Phi_0,\qquad J\qt=h\Phi_1,
 \qquad \beta=hb_0,\qquad \bt=hb_1.
\]
Thus the projected vectors have the common factor $h$,
and
\[
 D_{f,g}=h(g\Phi_0-f\Phi_1),
 \qquad
 \Psi_{f,g}=h(b_0|g|^2-b_1|f|^2).
\]
Multiplication by a nonzero scalar polynomial preserves
the coefficient kernel. The scalar sign polynomial
in~\eqref{eq:sign-test} has the additional factor $h^3$.
Since $h>0$ on the comparison interval, the scalar
sign conditions are unchanged.

The degree bounds for the minors and the scalar
sign polynomial follow from
$\deg\Phi_0,\deg\Phi_1\leq\ell-1$
and $\deg b_0,\deg b_1\leq\ell+1$.
Applying the coprime-factor argument in the proof of
Theorem~\ref{thm:subspace-classification}
to $Z^*\Phi_1$ and $Z^*\Phi_0$ gives
$\deg f,\deg g\leq\ell-1$.
\end{remark}

The relevant spectral projections also identify the directions orthogonal to both resolvent-vector families. Such directions can be added to any admissible subspace without changing the compressed quadratic form.

\begin{samepage}
\begin{corollary}\label{cor:invisible-subspace}
Let $\K=\spanop\{F^kq,F^k\qt:0\leq k<d\}$. Then $\Delta(x)$ is zero on $\K^\perp$ and has range in $\K$. If $\U$ is admissible, so is $\U+\K^\perp$.
\end{corollary}
\end{samepage}
\begin{proof}
The case $q=\zero$ is immediate. Otherwise use $\Lambda,\Phi_0,\Phi_1$ from Proposition~\ref{prop:reduced-degrees}; polynomial interpolation on the distinct eigenvalues of $F$ gives
\[
 \K=\spanop\{E_\lambda q,E_\lambda\qt:\lambda\in\Lambda\}.
\]
Both resolvent vectors and all coefficients of $g\Phi_0-f\Phi_1$ belong to $\K$. The assertions follow from~\eqref{eq:rank-two} and the coefficient description.
\end{proof}

\begin{remark}\label{rem:coefficient-kernels}
Every coefficient kernel in
Theorem~\ref{thm:subspace-classification}
contains $\K^\perp$.

Indeed, Lemma~\ref{lem:resolvent-span} shows that
the coefficients of $J(x)q$ and $J(x)\qt$ belong to $\K$.
Hence the coefficients of
\[
 D_{f,g}(x)=g(x)J(x)q-f(x)J(x)\qt
          =\sum_k d_kx^k
\]
also belong to $\K$.
For every $z\in\K^\perp$, we therefore have $d_k^*z=0$
for all $k$, and thus $\mathsf D_{f,g}^*z=\zero$.
Consequently,
\[
 \K^\perp\subseteq\ker\mathsf D_{f,g}^*.
\]
\end{remark}

If $q=\zero$, then $\Delta=\zero$ and every coupling subspace is admissible. In Example~\ref{ex:two-edges} with the additional diagonal weight $5$, the distinct eigenvalues of $F$ are $-1,0,1,5$ and $\Lambda=\{-1,0,1\}$. Thus $\ell=3$ and $d=4$.

\subsection{Maximum coupling dimension for a two-block matrix}
The coefficient kernels describe the subspaces associated with each rational relation. For the following matrix, the resolvent equations give an admissible subspace and a dimension bound directly. We state the result for a pair of bordered matrices with a common block.

Let $T\in\R^{r\times r}$ and $S\in\R^{s\times s}$ be symmetric and entrywise nonnegative, $a\in\R^r$ and $c\in\R^s$ be nonzero entrywise nonnegative vectors, and $b\in\R^s$ be entrywise nonnegative. For $\gamma>0$ and $\delta\in\R$, let $d=r+s$ and
\begin{equation}\label{eq:connected-bordered-pair}
 \begin{gathered}
 F=\begin{bmatrix}T&ab^T\\ba^T&S\end{bmatrix},\qquad
 q=\colvec{-\gamma a}{c},\qquad \qt=\colvec{\gamma a}{c},\\
 H_-=\begin{bmatrix}F&q\\q^T&\delta\end{bmatrix},\quad
 H_+=\begin{bmatrix}F&\qt\\\qt^T&\delta\end{bmatrix}.
 \end{gathered}
\end{equation}
For this pair, call $\Ran Z\subseteq\C^d$ admissible if the extensions of $H_-$ and $H_+$ by the same coupling $\colvec Z{\zero}$ satisfy the largest-eigenvalue comparison for every common Hermitian complementary block. This agrees with the preceding definition when the pair comes from a Kelmans transformation. Define the cyclic subspaces
\[
 \K_T=\spanop\{T^ja:0\leq j<r\},\qquad
 \K_c=\spanop\{S^jc:0\leq j<s\}.
\]
Denote their dimensions by $k_T=\dim\K_T$ and $k_c=\dim\K_c$.

The subspace $\C^r\oplus\K_c^\perp$ is admissible, and $k_T$ bounds the codimension of admissible subspaces that extend beyond it.

\begin{samepage}
\begin{theorem}\label{thm:rank-one-connection}
For the pair~\eqref{eq:connected-bordered-pair}, the subspace
$\U_0=\C^r\oplus\K_c^\perp$ is admissible. Every admissible $\U\subseteq\C^d$ satisfies
\[
 \U\subseteq\U_0\qquad\text{or}\qquad
 \operatorname{codim}_{\C^d}\U\geq k_T.
\]
Consequently, if $k_T\geq k_c$, the maximum admissible dimension is $d-k_c$. If $k_T>k_c$, then $\U_0$ is the unique admissible subspace of that dimension.
\end{theorem}
\end{samepage}
\begin{proof}
Let $h=\lambda_{\max}(H_+)$. Choose $t\geq0$ so that $H_-+tI$ and $H_++tI$ are positive semidefinite. Since $F\geew\zero$ and $\qt=|q|$, we have $|H_-+tI|\leew H_++tI$. For every real unit vector $y$, the Rayleigh quotients satisfy
\[
 y^T(H_-+tI)y\leq |y|^T(H_++tI)|y|\leq h+t.
\]
Taking the maximum over $y$ gives $\lambda_{\max}(H_-)\leq h$. For $x>h$, let $M=(xI-F)^{-1}$ and
\[
 s_-=x-\delta-q^TMq,\qquad s_+=x-\delta-\qt^TM\qt.
\]
Then $M\succ\zero$, $M\geew\zero$, and $0<s_+\leq s_-$. Block inversion and Theorem~\ref{thm:universal} show that $\Ran Z$ is admissible exactly when
\begin{equation}\label{eq:connected-compression}
 Z^*\left(\frac{1}{s_+}M\qt\qt^TM-\frac{1}{s_-}Mqq^TM\right)Z\succeq\zero
 \qquad(x>h).
\end{equation}
For the rest of the proof, let $x>h$ and let
\begin{align*}
 u&=(xI-T)^{-1}a,&v&=(xI-S)^{-1}c,&w&=(xI-S)^{-1}b,\\
 \alpha&=a^Tu,&\beta&=b^Tv,&\chi&=b^Tw,\qquad \varepsilon=1-\alpha\chi.
\end{align*}
Since $x>h\geq\rho(F)\geq\max\{\rho(T),\rho(S)\}$, both smaller resolvents are positive definite and entrywise nonnegative. Hence $\alpha>0$ and $\beta,\chi\geq0$. The Schur complement of $xI-T$ in $xI-F$ is $xI-S-\alpha bb^T\succ\zero$. Congruence by $(xI-S)^{-1/2}$ gives $I-\alpha\,w_0w_0^T\succ\zero$, where $w_0=(xI-S)^{-1/2}b$; if $b\ne\zero$, its eigenvalue in the direction $w_0$ is $1-\alpha\chi>0$. For $b=\zero$, $\varepsilon=1$. Thus $\varepsilon>0$ in both cases. For $R=M\colvec a{\zero}=\colvec{R_T}{R_S}$, the block equations are
\[
 (xI-T)R_T=a(1+b^TR_S),\qquad (xI-S)R_S=b(a^TR_T).
\]
Substitution and the identities $a^Tu=\alpha$, $b^Tw=\chi$ give
\[
 R_S=w(a^TR_T),\quad R_T=u(1+b^TR_S)=\frac{1}{\varepsilon}u,
 \quad R_S=\frac{\alpha}{\varepsilon}w.
\]
Also $(xI-F)\colvec{\zero}v=\colvec{-\beta a}{c}$. Consequently,
\[
 \begin{gathered}
 R=M\colvec a{\zero}=\frac{1}{\varepsilon}\colvec u{\alpha w},\qquad V=\colvec{\zero}v,\\
 Mq=(\beta-\gamma)R+V,\qquad
 M\qt=(\beta+\gamma)R+V.
 \end{gathered}
\]
If $Z_0$ spans $\U_0$, then $Z_0^*V=\zero$, since $v\in\K_c$. Thus
\[
 Z_0^*Mq=\frac{\beta-\gamma}{\beta+\gamma}Z_0^*M\qt.
\]
The scalar factor has modulus at most one. Together with $s_+\leq s_-$, Lemma~\ref{lem:rank-one-order} proves~\eqref{eq:connected-compression} for $Z_0$.

Now let $Z$ span any admissible $\U$, and let $W=\U^\perp$, $y=Z^*R$, and $z=Z^*V$. By~\eqref{eq:connected-compression} and Lemma~\ref{lem:rank-one-order}, the two projected border vectors are linearly dependent. They are $(\beta-\gamma)y+z$ and $(\beta+\gamma)y+z$. The coefficient matrix $\left[\begin{smallmatrix}\beta-\gamma&1\\\beta+\gamma&1\end{smallmatrix}\right]$ has determinant $-2\gamma\ne0$, so their dependence implies that $y,z$ are dependent for every $x>h$. If $z\equiv\zero$, Lemma~\ref{lem:resolvent-span} implies that $\U$ is orthogonal to $\{\zero\}\oplus\K_c$, so $\U\subseteq\U_0$.

Otherwise choose an entry $z_j$ that is not identically zero. Except at its finitely many zeros,
\[
 \omega(x)=z_j(x)R(x)-y_j(x)V(x)\in W,\qquad
 \operatorname{proj}_{\C^r}\omega(x)=\frac{z_j(x)}{\varepsilon(x)}u(x).
\]
Choose an open subinterval avoiding the zeros of $z_j$. Lemma~\ref{lem:resolvent-span} shows that the vectors $u(x)$ there span $\K_T$. Since $z_j(x)/\varepsilon(x)$ is nonzero, the projection of $W$ onto $\C^r$ contains $\K_T$. Its dimension is at most $\dim W$, so $\dim W\geq k_T$. This proves the dimension dichotomy, including complex $\U$. If $k_T\geq k_c$, both alternatives give $\dim\U\leq d-k_c$, attained by $\U_0$. If $k_T>k_c$, an admissible subspace of dimension $d-k_c$ must lie in $\U_0$ and hence equal it.
\end{proof}

The above theorem gives the maximum dimension when $k_T\geq k_c$, with uniqueness under the strict inequality. It applies to reducible $T,S$ and includes $b=\zero$. When $k_T<k_c$, the dichotomy gives the bounds
\[
 d-k_c\ \leq\ \max\{\dim\U:\U\text{ is admissible}\}\ \leq\ d-k_T.
\]

\subsection{A nonconstant proportionality factor}
\begin{example}\label{ex:nonconstant}
Let $G$ have vertices $0,1,2,3,4$ and edges $03,04,12$. Apply the graph Kelmans transformation toward $u=1$ from $v=3$. The new graph has edges $01,04,12$, so it is $P_4\cup K_1$ and $\rho(\widetilde A)=\varphi=(1+\sqrt5)/2$. The columns of $P$ are $(e_+,e_0,e_2,e_4)$, where the subscripts $0,2,4$ are original graph-vertex labels, and
\[
 F=\begin{bmatrix}
 0&1/\sqrt2&1/\sqrt2&0\\
 1/\sqrt2&0&0&1\\
 1/\sqrt2&0&0&0\\
 0&1&0&0
 \end{bmatrix},\qquad q=[0,-1/\sqrt2,1/\sqrt2,0]^T,\qquad\delta=0.
\]
Let the columns of $Z$ be the second and fourth standard basis vectors of $\C^4$. The columns of $PZ$ are then $e_0,e_4$, indexed by the original graph vertices. With $D(x)=2x^4-4x^2+1$, direct inversion yields
\[
 Z^*M\qt=\sqrt2\,D(x)^{-1}\colvec{x^3}{x^2},\qquad
 Z^*Mq=\frac{1-x^2}{x^2}Z^*M\qt.
\]
For $x>\varphi$, the factor $(1-x^2)/x^2$ has modulus less than one. Hence this two-dimensional coupling subspace is admissible. The displayed vector $Z^*M\qt$ is nonzero throughout that interval, so the proportionality factor is uniquely determined and depends on $x$.

Let $f(x)=1-x^2$ and $g(x)=x^2$.
A direct calculation gives
\[
 D_{f,g}(x)=\left[x^4-2x^2+\tfrac12,\ 0,\
 \frac{2x^5-4x^3+x}{\sqrt2},\ 0\right]^T.
\]
Its coefficient matrix has rank two, and the kernel of its adjoint is exactly the chosen subspace. Moreover,
\begin{align*}
 \beta&=x(x-1)(x+1)(x^2-2),\\
 \bt&=x(x^2-x-1)(x^2+x-1),\\
 \Psi_{f,g}&=x(x-1)(x+1)(2x^4-4x^2+1)>0\qquad(x>\varphi).
\end{align*}
Thus Theorem~\ref{thm:subspace-classification} includes this subspace through a rational factor of degree two.

It also attains the dimension bound in Corollary~\ref{cor:dimension-obstruction}. Here $\mathcal S$ is spanned by the second and third standard basis vectors of $\C^4$. The images of those two vectors under $F$ are
\[
 \left[\frac{1}{\sqrt2},0,0,1\right]^T,
 \qquad \left[\frac{1}{\sqrt2},0,0,0\right]^T.
\]
They lie in $\mathcal S^\perp$ and are independent. Corollary~\ref{cor:dimension-obstruction} gives dimension at most two, attained by the displayed subspace.
\end{example}
\section{Graph joins and Hermitian adjacency matrices}\label{sec:applications}

\subsection{Rank-one joins}
Let $A_1,\ldots,A_k$ be symmetric nonnegative matrices, let $\tau_i\geq0$, and choose nonnegative vectors $a_i$ of the corresponding orders. Embed each $a_i$ as a column supported on its own block, forming a matrix $U$. For a Hermitian matrix $H=[h_{ij}]_{i,j=1}^k$ with zero diagonal, consider
\begin{equation}\label{eq:join}
 \mathcal J=\diag(\tau_1A_1,\ldots,\tau_kA_k)+UHU^*.
\end{equation}
Thus the $(i,j)$ off-diagonal block is $h_{ij}a_i a_j^T$.

Theorem~\ref{thm:rank-one-monotonicity} can therefore be applied one diagonal block at a time, giving the following comparison.

\begin{samepage}
\begin{corollary}\label{cor:joins}
If $(a_i)_u=(a_i)_v$, applying a Kelmans transformation to $A_i$ toward $u$ from $v$ does not decrease $\lambda_{\max}(\mathcal J)$. The same holds for any sequence of such transformations.
\end{corollary}
\end{samepage}
\begin{proof}
For a single block, the assertion follows from Lemma~\ref{lem:inverse-comparison} and nonnegative scaling. Otherwise place the $i$th block first. Its coupling to the remaining blocks has the form $a_i b^*$ for a suitable complex vector $b$, and the remaining principal submatrix is a Hermitian matrix. Apply Theorem~\ref{thm:rank-one-monotonicity} at each step.
\end{proof}

When all $a_i=\one$, this includes the usual adjacency-matrix $H$-join. Repeated Kelmans transformations replace each graph block by a threshold graph of the same order and size~\cite{Csikvari11}, while weakly increasing the largest eigenvalue of~\eqref{eq:join}. The resulting extremal problem is therefore over threshold blocks of the prescribed orders and sizes.

\begin{remark}
The comparison concerns the largest eigenvalue. Closed-walk counts can decrease: a proper transformation sends $C_4$ to a triangle with a pendant edge, and $\operatorname{tr}(A^4)$ decreases from $32$ to $28$~\cite{Csikvari11}.
\end{remark}

The same argument applies if a single fixed diagonal unitary similarity puts both matrices into~\eqref{eq:join}, preserving the diagonal block structure and the rank-one factorization of every off-diagonal block with the same nonnegative vectors $a_i$. Eigenvalue invariance under this common similarity then gives the same comparison for the original pair.

\subsection{Higher-rank attachments}
Let $L$ be a simple graph on $r$ vertices, with adjacency matrix $T$, and let $\eta\in\{0,1\}^r$ be nonzero. Add vertices $u,v,w$, edges $uv,uw$, and edges $vj$ whenever $\eta_j=1$, with no other edges incident to $u,v,w$. Denote the resulting adjacency matrix by $A$. Let $\widetilde A$ be obtained by moving every edge $vj$ to $uj$, retaining $uv,uw$. This is a Kelmans transformation toward $u$ from $v$.

In the basis used below, this construction has the form of Theorem~\ref{thm:rank-one-connection}, with a cyclic subspace of dimension two in the second block. The theorem gives the following comparison and dimension bounds.

\begin{samepage}
\begin{corollary}\label{cor:higher-rank-attachments}
Let $m\geq1$, let $B\in\C^{(r+3)\times m}$ satisfy $B_{uj}=B_{vj}=B_{wj}=0$ for $1\leq j\leq m$, and let $C\in\C^{m\times m}$ be a Hermitian matrix. Then
\[
 \lambda_{\max}\begin{bmatrix}\widetilde A&B\\B^*&C\end{bmatrix}
 \geq
 \lambda_{\max}\begin{bmatrix}A&B\\B^*&C\end{bmatrix}.
\]
Let $k=\dim\spanop\{\eta,T\eta,\ldots,T^{r-1}\eta\}$. Among coupling subspaces contained in $(e_u-e_v)^\perp$, the largest admissible dimension is $r$ if $k\geq2$. If $k\geq3$, the subspace spanned by $\{e_j:j\in V(L)\}$ is the unique admissible subspace of dimension $r$.
\end{corollary}
\end{samepage}
\begin{proof}
In the orthonormal basis consisting of the vectors $e_j$, $j\in V(L)$, followed by $e_+,e_w,e_-$, the pair has the form~\eqref{eq:connected-bordered-pair} with $\gamma=1$, $\delta=-1$, and
\[
 a=\frac{1}{\sqrt2}\eta,\qquad b=\colvec10,\qquad
 S=\begin{bmatrix}1&1/\sqrt2\\1/\sqrt2&0\end{bmatrix},\qquad
 c=\colvec0{1/\sqrt2}.
\]
The edge $uv$ gives $e_+^TAe_+=1$ and $e_-^TAe_-=-1$. For $j\in V(L)$, the edges $vj$ give $e_j^TAe_+=\eta_j/\sqrt2$ and $e_j^TAe_-=-\eta_j/\sqrt2$; after the transformation, $e_j^T\widetilde A e_-=\eta_j/\sqrt2$. The edge $uw$ gives $e_+^TAe_w=e_-^TAe_w=1/\sqrt2$, unchanged by the transformation. The vectors $c,Sc$ are independent, so $k_c=2$ and $\U_0$ is precisely the subspace spanned by $\{e_j:j\in V(L)\}$. Apply Theorem~\ref{thm:rank-one-connection}, noting that $k_T=k$ and that a zero row corresponding to $e_-$ in this basis is equivalent to equal coupling rows at $u,v$ in the original basis.
\end{proof}

For $L=P_r$ and $\eta=e_1$, the vectors $T^j\eta$, $0\leq j<r$, are independent: $T^j\eta$ has a nonzero entry at distance $j$ from the endpoint and is zero at every greater distance. Thus this connected family has maximum admissible dimension $r$ for $r\geq2$, uniquely attained for $r\geq3$ among coupling spaces contained in $(e_u-e_v)^\perp$. In this family, the space $\K$ in Corollary~\ref{cor:invisible-subspace} is all of $\C^{r+2}$. Indeed, $\K=\spanop\{F^jq,F^j\qt:0\leq j<r+2\}$ contains $\colvec a{\zero}$ and $\colvec{\zero}c$. Since $b^Tc=0$, it also contains $\colvec{\zero}{Sc}$, hence the whole right block; applying $F$ and subtracting right-block components then gives $\colvec{T^ja}{\zero}$ for every $j$. Since $c,Sc$ span $\C^2$ and the vectors $T^ja$ span $\C^r$, these vectors fill both blocks, so $\K=\C^{r+2}$ and $\K^\perp=\{\zero\}$.

For unsigned attachments, Kelmans monotonicity gives the inequality~\cite{Csikvari11}. Here the coupling on $V(L)$ can have arbitrary signed or complex entries and arbitrary rank, and $C$ can be any Hermitian matrix. The dimension and uniqueness statements concern coupling spaces contained in $(e_u-e_v)^\perp$.

\subsection{A two-parameter Hermitian adjacency matrix}
Let $D$ be a loopless digraph, possibly with pairs of opposite arcs. For $\tau\in\R$ and $\zeta\in\C$ with $|\zeta|=1$, define the Hermitian matrix $H_{\tau,\zeta}(D)$ by zero diagonal and
\[
 (H_{\tau,\zeta}(D))_{ij}=
 \begin{cases}
 \tau,&i\to j\text{ and }j\to i,\\
 \zeta,&i\to j\text{ but not }j\to i,\\
 \overline\zeta,&j\to i\text{ but not }i\to j,\\
 0,&\text{neither arc is present}.
 \end{cases}
\]
Here the $(i,j)$ entry corresponds to the arc direction $i\to j$ in the second case. Several standard Hermitian adjacency conventions fit this framework~\cite{GuoMohar,LiuLi,Mohar}.

Let $\overline D$ denote the arc complement: each ordered pair of distinct vertices is an arc of $\overline D$ exactly when it is not an arc of $D$. Complementation interchanges digons with nonedges and reverses each single arc. Consequently,
\begin{equation}\label{eq:arc-complement}
 H_{\tau,\zeta}(\overline D)=\tau(\one\one^T-I)-H_{\tau,\zeta}(D)
 \quad\text{for every }D
 \quad\Longleftrightarrow\quad\tau=2\Re\zeta.
\end{equation}
Both sides have zero diagonal. At a digon or nonedge, their entries agree for every $\tau$. At a single arc, agreement is equivalent to $\overline\zeta=\tau-\zeta$, or $\tau=2\Re\zeta$, proving both directions of~\eqref{eq:arc-complement}.

\begin{definition}\label{def:module}
A vertex subset $S$ is a \emph{symmetric module} if every pair inside $S$ is either a digon or a nonedge, and all vertices in $S$ have the same in-neighbors and the same out-neighbors outside $S$.
\end{definition}

Let $G_S$ be the undirected graph on $S$ whose edges are its digons. Ordering $S$ first gives
\[
 H_{\tau,\zeta}(D)=\begin{bmatrix}\tau A(G_S)&\one b^*\\b\one^T&C\end{bmatrix}
\]
for a complex vector $b$ and a Hermitian matrix $C$. Thus the module condition gives the rank-one coupling in Section~\ref{sec:rank-one}.

\begin{samepage}
\begin{theorem}\label{thm:digraph}
Suppose $S$ is a symmetric module. Apply a graph Kelmans transformation to $G_S$, and let $\widetilde D$ be obtained by moving the corresponding digons in $D$, retaining all other arcs. If $\tau\geq0$, then
\begin{align*}
 \lambda_{\max}(H_{\tau,\zeta}(\widetilde D))&\geq\lambda_{\max}(H_{\tau,\zeta}(D)),\\
 \lambda_{\max}(H_{\tau,\zeta}(\overline {\widetilde D}))&\geq\lambda_{\max}(H_{\tau,\zeta}(\overline D)).
\end{align*}
For each $\zeta$ on the unit circle, the two inequalities hold for every such digraph and transformation exactly when $\tau\geq0$.
\end{theorem}
\end{samepage}
\begin{proof}
If the exterior of $S$ is nonempty, apply Theorem~\ref{thm:rank-one-monotonicity} to the displayed module block form; if it is empty, use Lemma~\ref{lem:inverse-comparison} and scaling. In the arc complement, $S$ is again a symmetric module and its internal graph is $\overline{G_S}$. Complementation turns a Kelmans transformation toward $u$ from $v$ into a Kelmans transformation toward $v$ from $u$~\cite{Csikvari11}. Applying the same argument to this reverse transformation proves the second inequality for every $|\zeta|=1$ and $\tau\geq0$.

For $\tau=-c<0$, take $S$ to be all vertices and replace every edge of the path in Theorem~\ref{thm:rank-one-monotonicity} by a digon. Its transformation gives a digon triangle and an isolated vertex. The largest eigenvalue decreases from $c\varphi$ to $c$, which completes the proof.
\end{proof}

Mohar's matrix $N_\alpha$ assigns $\alpha$ to the $(i,j)$ entry and $\overline\alpha$ to the $(j,i)$ entry of a single arc $i\to j$, and $\alpha+\overline\alpha$ to a digon, with $|\alpha|=1$ and $\Re\alpha\geq0$~\cite{Mohar}. Thus
\[
 N_\alpha(D)=H_{2\Re\alpha,\alpha}(D).
\]
Theorem~\ref{thm:digraph} applies throughout this parameter range. If the same formula is extended to the entire unit circle, the universal module comparison holds exactly when $\Re\alpha\geq0$.
\section{AI Declaration}\label{sec:AIdeclaration}

ChatGPT (GPT) assisted with the exploration of proof strategies and the formulation and testing of intermediate claims. The authors directed the research, selected the methods, organized the arguments, independently checked all proofs, and take responsibility for the mathematical content and presentation.

The work with GPT began in early August 2026 with spectral comparison under Kelmans transformations in Hermitian block matrices. An orthogonal change of basis identified the common principal block before and after the transformation. Schur complements~\cite{HornZhang,LiMathias} then reduced comparison for every common complementary block to the order of compressed resolvents. This led to the study of coupling subspaces. Range inclusion and cyclic subspaces~\cite{Douglas,Godsil} were used to describe these subspaces and their dimensions, followed by applications to graph joins, attachments, and Hermitian adjacency matrices.

\end{document}